\documentclass[12pt]{article}
\usepackage{amsmath, amssymb}
\usepackage{amssymb,pb-diagram}
\usepackage{amscd}

\usepackage{fancybox}

\newcommand{\ZZ}{\mathbb Z}
\newcommand{\PP}{\mathbb P}
\newcommand{\QQ}{\mathbb Q}
\newcommand{\CC}{\mathbb C}

\newcommand{\MW}{\mathop {\rm MW}\nolimits}
\newcommand{\Gal}{\mathop {\rm Gal}\nolimits}

\newcommand{\NS}{\mathop {\rm NS}\nolimits}

\newcommand{\Red}{\mathop {\rm Red}\nolimits}

\newcommand{\Pic}{\mathop {\rm Pic}\nolimits}

\newcommand{\Supp}{\mathop {\rm Supp}\nolimits}

\def\miya#1{\mathrel{\mathop{\rightarrow}\limits^#1}}

\newcommand{\id}{\mathop {{\rm id}}\nolimits}

\newtheorem{thm}{Theorem}[section]
\newtheorem{cor}{Corollary}[section]
\newtheorem{prop}{Proposition}[section]
\newtheorem{lem}{Lemma}[section]
\newtheorem{defin}{Definition}[section]
\newtheorem{exmple}{Example}[section]
\newtheorem{rem}{Remark}[section]
\newtheorem{qz}{Question}[section]
\newtheorem{prbm}{Problem}[section]
\newcommand{\qed}{\hfill $\Box$}

\newcommand{\proof}{\noindent{\textsl {Proof}.}\hskip 3pt}
\newcommand{\proofend}{\qed \par\smallskip\noindent}

\renewcommand{\thesubparagraph}{\theparagraph.\@arabic\c@subparagraph}
\begin{document}
  
  \title{ \bf Splitting curves on a rational ruled surface,  the Mordell-Weil groups of
  hyperelliptic fibrations and
  Zariski pairs}

\author{ Hiro-o TOKUNAGA}

\date{}


\maketitle

\begin{abstract}  
 Let $\Sigma$ be a smooth projective surface,  let $f' : S' \to \Sigma$ be a double 
 cover of $\Sigma$ and let $\mu : S \to S'$ be the canonical resolution.  Put $f = f'\circ\mu$.
 An irreducible curve $C$ on $\Sigma$ is said to be a 
 splitting curve with respect to $f$ if $f^*C$ is of the form  $C^+ + C^- + E$, where
 $C^- = \sigma_f^*C^+$,  $\sigma_f$ being  the covering transformation of $f$ and 
 all irreducible components of $E$ are contained in the exceptional set of $\mu$. 
 In this article, we show 
 that a kind of  ``reciprocity"  of splitting curves holds for a certain pair of  curves on rational ruled surfaces. As an application, we consider the topology of the complements of certain curves
 on rational ruled surfaces.
 
\end{abstract}

  \Large
  
  {\bf Introduction}

  \normalsize

  Let $Y$ be a smooth projective variety, and let $X'$ be a double cover, i.e., a normal variety
  with a finite surjective morphism $f' : X' \to Y$ of $\deg f' = 2$. We denote its smooth model by
  $X$ and its resolution
  by $\mu_{f'} : X \to X'$ and  put $f = f'\circ\mu_{f'}$. Let $\sigma_{f'}$ be
 the covering transformation of $f'$ and we assume that $\sigma_{f'}$ inducdes an involution
  $\sigma_f$ on $X$. 
  \begin{defin}\label{def:q-residue-div}{\rm Let ${\mathcal D}$ be an irreducible divisor on $Y$.
  \begin{enumerate}
 \item[(i)] We say that $D$ is a splitting divisor with respect to $f : X \to Y$ if
  $f^*{\mathcal D}$ is of the form
\[
f^*{\mathcal D} = {\mathcal D}^+ + {\mathcal D}^- + E,
\]
where ${\mathcal D}^- = \sigma_f^*{\mathcal D}^+$, and irreducible components of $E$ are all exceptional divisors
of $\mu_{f'}$. In case of $\dim Y =2$, we call ${\mathcal D}$ a splitting curve.

\item[(ii)] If $f'$ is uniquely determined by the branch locus
 $\Delta_{f'}$ of $f'$,  then we say that $D$ is a splitting divisor
with respect to $\Delta_{f'}$ (see \S 1 for the terminologies for branched covers).
Note that if $Y$ is simply connected, then any double cover is uniquely determined by its branch locus.
\end{enumerate}
}
\end{defin}

\begin{rem}{\rm In the case of $\dim Y = 2$, if ${\mathcal D}$ is a splitting curve on $Y$ with respect to
$f : X \to Y$, ${\mathcal D}$ does not meet $\Delta_{f'}$ transversely at smooth parts of 
${\mathcal D}$ and $\Delta_{f'}$.  
}
\end{rem}

In this article, we study splitting curves on a smooth surface from two different viewpoint:
 
 \begin{enumerate}
 
 \item[(I)] The study of dihedral covers and their application.
 
 \item[(II)] A problem motivated by elementary number theory, i.e., to formulate ``reciprocity
 law" of double covers.
 
 \end{enumerate}
 
In order  to illustrate how the notion of splitting divisors works in the study of dihedral covers,
let us recall some results in \cite{artal-tokunaga} (For
notations and terminology, see \S 1):

\medskip

Let $C$ be a smooth conic in $\PP^2$ and let $f : Z \to \PP^2$ be the double cover
of $\PP^2$ with branch locus $\Delta_f = C$. Note that $Z \cong \PP^1\times \PP^1$.
We denote the class of a divisor in $\Pic(Z)$ by a pair of integers $(a, b)$. Then we have:

\begin{prop}\label{prop:at}{{\rm (\cite[Proposition 2]{artal-tokunaga} )}
For  an irreducible curve $D$ in $\PP^2$, there exists a $D_{2n}$-cover $\pi : X \to \PP^2$
such that 
\begin{itemize}

 \item $\Delta_{\pi} = C + {\mathcal D}$, and
 \item the ramification index along $C$ (resp. ${\mathcal D}$)  is $2$ (resp. $n$) (we say that $f$ is branched at 
 $2C+ n{\mathcal D}$ if these conditions on the ramification are satisfied),
 
 \end{itemize}
 
  if and only if both of  two condititons below are stasified:
  
  \begin{enumerate}
   \item[(i)] ${\mathcal D}$ is a splitting curve with respect to $C$.
   
   \item[(ii)] If we put $f^*{\mathcal D} = {\mathcal D}^+ + {\mathcal D}^-$ and  denote the class of ${\mathcal D}^+$ by $(a, b)$, then
   $a - b$ is divisible by $n$.
  
   \end{enumerate}
   }
   \end{prop}
   
     In \cite{artal-tokunaga}, we construct nodal rational curves ${\mathcal D}_1, \ldots, {\mathcal D}_k$ of degree $m$ in 
   $\PP^2$ as follows:
   
   \begin{itemize}
   \item $k$ is an integer not exceeding $m/2$.
   
   \item For each $i$, ${\mathcal D}_i$ is tangent to $C$ at $m$  distinct smooth points of 
   ${\mathcal D}_i$.
   
   \item $f^*{\mathcal D}_i = {\mathcal D}_i^+ + {\mathcal D}_i^-, {\mathcal D}^+_i \sim (m-i, i),
    {\mathcal D}^-_i \sim (i, m-i)$.
   
   \end{itemize}
   
  By Propositin~\ref{prop:at}, we can show that the following statement:

   \begin{prop}\label{prop:at-kplet}{
    For any $i, j (i \neq j)$, there exists no homeomorphism $f : \PP^2 \to \PP^2$ such that 
   $f(C\cup {\mathcal D}_i) = C\cup {\mathcal D}_j$ for any $i$ and $j$, i.e.,
there is no homeomorphism between pairs $(\PP^2, C\cup {\mathcal D}_i)$ and $(\PP^2, C\cup 
{\mathcal D}_j)$.
Namely $(C\cup {\mathcal D}_1, \ldots, C\cup {\mathcal D}_k)$ is a Zarski $k$-plet (See \cite{act} for the definition
of a Zariski $k$-plet or a Zariski pair).    
 }
 \end{prop}  
     



%



In Proposition~\ref{prop:at-kplet}, one of clues  is that ${\mathcal D}_i$ is rational, and it is a splitting curve with respect
to $C$.  In this article, we consider the case that ${\mathcal D}_i$ is {\it non-rational}.
\medskip

We also remark that in \cite{shimada} Shimada intensively studied splitting curves of degree $\le 2$
on $\PP^2$
with respect to sextic curves with simple singularities.  Such splitting curves play
essential role to classify so called \textit{lattice Zariski pairs}.

\medskip

As for the viewpoint (II), let us recall a
fact from number theory:

Let $m$ be a square free integers and put $K = \QQ(\sqrt {m})$. We denote the ring
of integers of $K$ by $O_K$ and the discriminant of $K$ by $\delta_K$. Let $p$ be 
an odd prime with $p \not| \delta_K$ and let $(p)$ be  the ideal  of $O_K$ generated
by $p$. Then the statements below hold (See \cite[Proposition 13.1.3]{ireland-rosen}, p.190, for example)
:

\begin{enumerate}

\item[(i)] If $x^2 \equiv m \bmod p$ is solvable in $\ZZ$, then $(p) = {\mathfrak p}_1{\mathfrak p}_2$, where ${\mathfrak p}_1$ and ${\mathfrak p}_2$ are distinct prime ideals in $O_K$.

\item[(ii)] If $x^2 \equiv m \bmod p$ is not solvable  in $\ZZ$, then $(p)$ is a prime ideal in 
$O_K$.

\end{enumerate}

Hence  whether $(p)$ splits or not depends on the solvability of $x^2 \equiv m \bmod p$. 
Moreover law of quadratic reciprocity gives a relation between the solvability of 
$x^2 \equiv q \bmod p$ and that of $x^2 \equiv p \bmod q$ for odd primes.

These facts suggest us to formulate the following problem:

\begin{prbm}\label{pbm:basic}{\rm Let $\Sigma$ be a smooth projective surface. Let 
$D_1, D_2$ and $D_3$ be reduced divisors on $\Sigma$, We denote the irreducible decompostion
of $D_i$ $(i = 1, 2)$ by $D_i = \sum_j D_{i, j}$ $(i = 1, 2)$, respectively. Suppose that there
exist double covers $p_i : {\mathcal S}_i \to \Sigma$ with $\Delta_{p_i} = D_i + D_3$ for $i = 1, 2$.
Is there any law to determine whether $p_1^*D_{2, j}$ splits or not in terms of some propeties of 
${\mathcal S}_2$?
}
\end{prbm}

   In this article, keeping the viewpoint (I), in particular, application to the study of  Zariski pairs (\S 5),
     in mind, we consider Problem~\ref{pbm:basic} under the 
   following setting:   
   \par\medskip
   
   Let  $\Sigma_{d}$ be the Hirzebruch surface of degree $d$. Throughout Introduction,
   we assume that \textbf{$d$ is even}.
  $\Delta_{0,d}$  denotes the negative section, i.e., the section whose self-intersection number is
   $-d$ and  $F_d$ denotes  a fiber of $\Sigma_{d} \to \PP^1$. Note that $\Pic(\Sigma_{d}) = \ZZ\Delta_{0,d}\oplus \ZZ F_d$.

    Let  $T_d$ be an irreducible divisor on $\Sigma_d$ such that 
    
    \begin{enumerate}
    \item[(i)] $T_d \sim (2g+1)(\Delta_{0,d} + d F_d)$ ($g \ge 1$)
   and
   
   \item[(ii)]  $T_d$ has only nodes  (resp. at worst simple singularities ) if $g \ge 2$ (resp.
   $g = 1$).
   
   \end{enumerate}

 Let $\Delta$ be  a section of $\Sigma_{d}$
 such that 
 
 \begin{enumerate}
 
\item[(i)] $\Delta \sim \Delta_{0,d} + d F_d$ and,
 
\item[(ii)] for all $x \in T_d\cap \Delta$, $x$ is a smooth 
 point of $T_d$ and the local intersection number $(\Delta\cap T_d)_x$ at $x$ is even.
 
 \end{enumerate}
 
 Let  $q_d : W_d \to \Sigma_{d}$ be a double cover branched at $2(\Delta_{0,d} + \Delta)$. 
 Let $p'_{d} : S'_{d} \to \Sigma_{d}$ be a double cover branched at 
 $2(\Delta_{0, d} + T_d)$ and 
 let $\mu_{d} : S_{d} \to S'_{d}$ be the canonical resolution in the sense of \cite{horikawa}.
 Put $p_{d} = p'_{d}\circ\mu_{d}$.  We can easily check the following properties:
 \begin{itemize}
  \item $W_{d}$ is the Hirzebruch surface of degree $d/2$.
  
  \item The composition $\varphi_{d} : S_{d} \to \Sigma_{d} \to \PP^1$ gives
   a hyperelliptic fibration of genus $g$ on $S_{d}$. The preimage of $\Delta_{0,d}$ in 
   $S_d$  gives a section
   $O$. We denote the Mordell-Weil group of the Jacobian of  the geneiric fiber $(S_{d})_{\eta}$
   by $\MW({\mathcal J}_{S_d})$, $O$ being the zero element.
   
   \item $p^*_{d}\Delta$ is of the form
   \[
   p^*_{d}\Delta = s^+ + s^-,
   \]
   where $s^{\pm}$ are sections of $\varphi_{d}$ and $s^- = - s^+$ in $\MW({\mathcal J}_{S_d})$.  In particular, $\Delta$ is a splitting curve with respect to $p_{d}$.
   
   \end{itemize}
   
   Now we are in position to state our main result in this article:
   
   \begin{thm}\label{thm:qr-1}{Under the notation as above, $T_d$ is a splitting curve with respect to $q_{d}$, if and only if $s^+$ is $2$-divisible in
   $\MW({\mathcal J}_{S_d})$.
   }
   \end{thm}

 Note that Theorem~\ref{thm:qr-1} can be regarded as an analogy of the ``reciprocity."  In fact,
 we may interpret it as a double cover version of 
 \[
 \left ( \frac 2p \right ) = (-1)^{\frac {p^2 -1}8}.
 \]
 
 It may
be interesting problem to consider the ``reciprocity" of double covers
 under more general setting. 

As an application of Theorem~\ref{thm:qr-1}, we study branched covers with
branch locus $\Delta_{0,d}+ \Delta + T_d$ and the topology of the pair $(\Sigma_d, \Delta_{0, d}\cup \Delta\cup T_d)$.
Here is our statement:

\begin{thm}\label{thm:application}{We keep the notation as before. There exists a $D_{2n}$-cover
branched at $2(\Delta_{0,d} + \Delta) + nT_d$ for any odd number $n$, if and only if $s^+$ is $2$-divisible in 
$\MW({\mathcal J}_{S_d})$
}
\end{thm}

An interesting corollary to Theorem~\ref{thm:application} is as follows:

\begin{cor}\label{cor:zpair}{
Let $(\Delta_1\cup T_{d,1}, \Delta_2\cup T_{d, 2})$ be a pair of  reduced divisors on $\Sigma_d$ such that both of $\Delta_i \cup T_{d,i}$ $( i = 1, 2)$  satisfy the conditions of $\Delta$ and $T_d$ in Theorem~\ref{thm:qr-1}. 
 Let $p_{d, i} : S_{d,i} \to \Sigma_d$ $(i = 1,2 )$ be the
canonical resolutions of 
the double covers with branch locus $\Delta_{0,d} + T_{d,i}$ $(i = 1, 2)$, respectively.
Let $s_i^+$ $(i = 1,2)$ be sections coming from $\Delta_{i}$ $(i = 1, 2)$, respectively.
If 
$(i)$ $s_1^+$ is $2$-divisible in $\MW({\mathcal J}_{S_{d,1}})$ and $(ii)$ $s_2^+$ is not $2$-divisible in
$\MW({\mathcal J}_{S_{d,2}})$, then there is no homeomorphism $h : \Sigma_d \to \Sigma_d$ such that
$h(\Delta_{0, d}) = \Delta_{0,d}$, 
$h(\Delta_1) = \Delta_2$ and $h(T_{d, 1}) = T_{d, 2}$.
}
\end{cor}

This article consists of 4 sections. In \S 1, we  summarize on Galois covers,
especially, dihedral covers, 
 and the Mordell-Weil group of the Jacobian
of the generic fiber of a fibered surface.  
We prove Theorems~\ref{thm:qr-1} and ~\ref{thm:application} in \S 2 and \S 3,
respectively. In \S 4, we give some 
explicit examples in the case when $d = 2, g=1$. In the last section,  we will see that the
study of splitting curves gives some examples of Zariski pairs, which is related our viewpoint I.

 

\section{Preliminaries}

\subsection{Summary on Galois  covers}

\textit{1. Generalities}
\medskip

Let $G$ be a finite group. Let $X$ and $Y$ be a normal projective varieties. We call $X$ a $G$-cover, if there exists a finite
surjective morphism $\pi : X \to Y$ such that the finite field extension given by
$\pi^* : \CC(Y) \hookrightarrow \CC(X)$ is a Galois extension with 
$\Gal(\CC(X)/\CC(Y)) \cong G$. We denote the branch locus of $\pi$ by
$\Delta_{\pi}$. $\Delta_{\pi}$ is a reduced divisor if $Y$ is smooth(\cite{zariski}).
 Let $B$ be a reduced divisor on $Y$ and we denote its irreducible decomposition by $B = \sum_{i=1}^r B_i$. We say that a $G$-cover 
$\pi : X \to Y$ is branched at $\sum_{i=1}^r e_iB_i$ if 

\begin{enumerate}
\item[(i)] $\Delta_{\pi} = \Supp (\sum_{i=1}^r e_i B_i)$ and
\item[(ii)] the ramification index along $B_i$ is $e_i$ for $1 \le i \le r$.
\end{enumerate}

\bigskip

\textit{2. Cyclic covers and double covers}

\medskip

Let $\ZZ/n\ZZ$ be a cyclic group of order $n$. We call a $\ZZ/n\ZZ$- (resp. a $\ZZ/2\ZZ$-) cover
by an $n$-cyclic (resp. a double) cover. 
We here summarize some facts
on cyclic and double covers. We first remark the following fact on cyclic covers.

\medskip

\textbf{Fact:} Let $Y$ be a smooth projective variety and $B$ a reduced divisor on $Y$. If
there exists a line bundle ${\mathcal L}$ on $Y$ such that 
$B \sim n{\mathcal L}$, then we can construct a hypersurface $X$ in the total space, $L$, of 
${\mathcal L}$ such that 

\begin{itemize}

\item $X$ is irreducible and normal, and

\item $\pi := \mbox{pr}|_X$ gives rise to an $n$-cyclic cover, pr being the canonical projection
$\mbox{pr} : L \to Y$.
\end{itemize}
(See \cite{bpv} for the above fact.)

\medskip

As we see in \cite{tokunaga90}, cyclic covers are not always realized as a hypersurface of the
total space of a certain line bundle. 
As for double covers, however, the following lemma holds.

\begin{lem}\label{lem:db-1}{Let $f : X \to Y$ be a double cover of a smooth projective
variety with $\Delta_f = B$, then there exists a line bundle ${\mathcal L}$ such that
$B \sim 2{\mathcal L}$ and $X$ is obtained as a hypersurface of the total space, $L$, of
${\mathcal L}$ as above.
}
\end{lem}

\proof  Let $\varphi$ be a rational function in $\CC(Y)$ such that 
$\CC(X) = \CC(Y)(\sqrt{\varphi})$. By our assumption,  the divisor of $\varphi$ is of the form
\[
(\varphi) = B + 2D,
\]
where $D$ is a divisor on $Y$. Choose ${\mathcal L}$ as the line bundle determined by $-D$. This implies our 
statement.
\proofend

\medskip

By Lemma~\ref{lem:db-1}, note that any double cover $X$ over $Y$ is determined by the pair
$(B, {\mathcal L})$ as above. 

\bigskip

\textit{3. Dihedral covers}

\medskip
 
 We next explain dihedral covers briefly.
 Let $D_{2n}$ be a dihedral group of order $2n$ given by
$\langle \sigma, \tau \mid \sigma^2 = \tau^n = (\sigma\tau)^2 = 1\rangle$.  In \cite{tokunaga94},
we developed a method of dealing with $D_{2n}$-covers. We need to introduce some notation
in order to describe it.

Let $\pi : X \to Y$ be a $D_{2n}$-cover. By its definition, $\CC(X)$ is a $D_{2n}$-extension
of $\CC(Y)$. Let $\CC(X)^{\tau}$ be the fixed field by $\tau$. We denote the $\CC(X)^{\tau}$-
normalization by $D(X/Y)$. We denote the induced morphisms by
$\beta_1(\pi) : D(X/Y) \to Y$ and $\beta_2(\pi) : X \to D(X/Y)$. Note that $X$ is a $\ZZ/n\ZZ$-cover 
of $D(X/Y)$ and $D(X/Y)$ is a double cover of $Y$ such that $\pi = \beta_1(\pi)\circ\beta_2(\pi)$:

\[
\begin{diagram}
\node{X}\arrow[2]{s,l}{\pi}\arrow{se,t}{\beta_2(\pi)}\\
 \node{}
\node[1]{D(X/Y)}\arrow{sw,r}{\beta_1(\pi)} \\ 
\node{Y}
\end{diagram}
\]

In \cite{tokunaga94}, we have the following results for $D_{2n}$-covers ($n$:odd).

\begin{prop}\label{prop:di-suf}{Let $n$ be an odd integer with $n \ge 3$. 
Let $f : Z \to Y$ be a double cover of a smooth projective variety $Y$, and assume that $Z$ is smooth. 
 Let $\sigma_f$ be the covering transformation of $f$. Suppose that  there exists a pair 
 $(D, {\mathcal L})$ of an effective divisor and a 
line bundle on $Y$ such that
\begin{enumerate}
\item[(i)] $D$ and $\sigma_f^*D$ have no common components,

\item[(ii)] if $D = \sum_ia_iD_i$ denotes the irreducible decomposition of $D$, then 
$0 < a_i \le (n-1)/2$ for every $i$;and the greatest common divisor of the $a_i$'s and $n$ is $1$, and

\item[(iii)] $D - \sigma_f^*D$ is linearly equivalent to $n{\mathcal L}$.
\end{enumerate}
Then there exists a $D_{2n}$-cover, $X$, of $Y$ such that (a) $D(X/Y) = Z$,
(b) $\Delta(X/Y) =\Delta_f \cup f(\Supp D)$ and (c) the ramification index along $D_i$ is
$n/\gcd(a_i, n)$ for $\forall i$. 
}
\end{prop}

For a proof, see \cite{tokunaga94}. A corollary related to a splitting divisor on $Y$, we have
the following:

\begin{cor}\label{cor:di-suf}{Let ${\mathcal D}$ be a splitting divisor on $Y$ with respect to $f :
Z \to Y$. Put $f^*{\mathcal D} = {\mathcal D}^+ + {\mathcal D}^-$. If there exists a line bundle
${\mathcal L}$ on $Z$ such that ${\mathcal D}^+ - {\mathcal D}^- \sim n{\mathcal L}$ for an
odd number $n$,
then there exists a $D_{2n}$-cover $\pi : X \to Y$ branched at $2\Delta_f + n{\mathcal D}$.
}
\end{cor}

Conversely we have the following:

\begin{prop}\label{prop:dinec}{Let $\pi : X \to Y$ be a $D_{2n}$-cover ($n \ge 3$, $n$: odd) of $Y$
and let $\sigma_{\beta_1(\pi)}$ be the involution on $D(X/Y)$ determined by
the covering transformation of $\beta_1(\pi)$. Suppose that
$D(X/Y)$ is smooth. Then there exists a pair of an effective divisor and a line bundle $(D, {\mathcal L})$
on
$D(X/Y)$ such that 

\begin{enumerate}
\item[(i)] $D$ and $\sigma_{\beta_1(\pi)}^*D$ have no common component,
\item[(ii)] if $D = \sum_ia_iD_i$ denotes the decomposition into irreducible components, then
$0 \le a_i \le (n-1)/2$ for every $i$,

\item[(iii)] $D - \sigma_{\beta_1(\pi)}^*D \sim n{\mathcal L}$, and

\item[(iv)] $\Delta_{\beta_2(\pi)} = \Supp(D+\sigma_{\beta_1(\pi)}^*D)$.
\end{enumerate}
}
\end{prop}

For a proof, see \cite{tokunaga94}.

\begin{cor}\label{cor:di-nec}{ Let $D_i$ be an arbitrary irreducible component of $D$ in 
Proposition~\ref{prop:dinec}. The image $\beta_1(\pi)(D_i)$ is a splitting divisor with respect
to $\beta_1(\pi) : D(X/Y) \to Y$
}
\end{cor}

\subsection{A review on the Mordell-Weil groups for fibrations over curves}

In this section, we review the results on the Mordell-Weil group and the Mordell-Weil
lattices studied by Shioda in \cite{shioda90, shioda99}.

Let $S$ be a smooth algebraic surface with  fibration $\varphi : S \to C$ of genus $g (\ge 1)$ curves over a smooth curve $C$. Throughout this section, we assume that

\begin{itemize}
 \item $\varphi$ has a section $O$ and
 \item $\varphi$ is relatively minimal, i.e., no $(-1)$ curve is contained
 in any fiber.
 \end{itemize}
 
 Let $S_{\eta}$ be the generic fiber of $\varphi$ and let $K =\CC(C)$ be the rational
 function field of $C$. $S_{\eta}$ is regarded as a curve of genus $g$ over $K$.
 
 Let ${\mathcal J}_S:= J(S_{\eta})$ be the Jacobian variety of $S_{\eta}$. We denote the set of rational points
 over $K$ by $\MW({\mathcal J}_S)$. By our assumption, $\MW({\mathcal J}_S) \neq \emptyset$ and it is 
 well-known that $\MW({\mathcal J}_S)$ has the structure of an abelian group.

 Let $\NS(S)$ be the N\'eron-Severi group of $S$ and  let $Tr(\varphi)$ be the subgroup
 of $\NS(S)$ generated by $O$ and irreducible components of fibers of $\varphi$. Under these
 notation, we have: 
 
 \begin{thm}\label{thm:shioda}{ If the irregularity of $S$ is equal to $C$, then
 $\MW({\mathcal J}_S)$ is a finitely generated abelian group such that
 \[
 \MW({\mathcal J}_S) \cong \NS(S)/Tr(\varphi).
 \]
 }
 \end{thm}
 
 See \cite{shioda90, shioda99} for a proof.
 
 \bigskip
 
 Let $p_d :S_d \to \Sigma_d$ be the double cover of $\Sigma_d$ with branch locus
 $\Delta_{0,d} + T_d$ as in Introduction. Then we have

 \begin{lem}\label{lem:non-etale}{Let $p_d : S_d \to \Sigma_d$ be the double cover as before.
 There exists no unramified cover of $S_d$. In particular, $\Pic(S_d)$ has no torsion element.
 }
 \end{lem}
 
 \proof  By Brieskorn's results on the simultaneous resolution of rational double points, we
 may assume that $T_d$ is smooth. Since the linear system $|T_d|$ is base point free,
 it is enough to prove our statement for one special case. Chose an affine open set $U_d$
 of $\Sigma_d$ isomorphic to $\CC^2$ with a coordinate $(x, t)$ so that a curve  $x = 0$ gives
 rise to a section linear equivalent to $\Delta_{\infty}$.
 Choose  $T_d$  whose defining equation in $U_d$ is
 \[
 T_d: F_{T_d} =  x^{2g+1} -  \Pi_{i=1}^{(2g+1)d}(t - \alpha_i) = 0,
 \]
 where $\alpha_i$ ($i = 1,\ldots, (2g+1)d$) are distinct complex numbers.  Note that
 \begin{itemize}
 
 \item  $T_d$ is smooth,
  
\item   singulair fibers of $\varphi$ are over $\alpha_i$ ($i = 1, \ldots, (2g+1)d$), and

\item all the singular fibers are irreducible rational curves with unique singularity isomorphic to
 $y^2 - x^{2g+1}= 0$. 
 
 \end{itemize}
  Suppose that $\widehat{S}_d \to S_d$ be any unramified cover, and let $\hat {g} : \widehat{S}_d
  \to \PP^1$ be the induced fiberation.  We claim that $\hat {g}$ has a connected fiber. Let
  $\widehat{S}_d \miya{{\rho_1}} C \miya{{\rho_2}} \PP^1$ be the Stein factorization and 
  let $\widehat {O}$ be a section coming from $O$. Then $\deg (\rho_2\circ\rho_1)|_{\hat {g}} = \deg \hat {g}|_{\widehat {O}} = 1$, and $\hat g$ has a connected fiber. 
  
  On the other hand, since all the singular fibers of $g$ are simply connected, all fibers over
  $\alpha_i$ ($i = 1, \ldots, (2g+1)d$) are disconnected. This leads us to a contradiction.
  \proofend

  \begin{cor}\label{cor:mw}{
 The irreguarity $h^1(S_d, {\mathcal O}_{S_d})$ of $S_d$ is $0$. In particular, 
 \[
 \MW({\mathcal J}_{S_d}) \cong \NS(S_d)/Tr(\varphi),
 \]
 where $Tr(\varphi)$ denotes the subgroup of $\NS(S_d)$ introduced as above.
 }
 \end{cor}
 
 \proof By Lemma~\ref{lem:non-etale}, we infer that $H^1(S_d, \ZZ) = \{0\}$. Hence
 the irregularity of $S_d$ is $0$.
 
 \proofend
 

\section{Proof of Theorem~\ref{thm:qr-1}}

Let us start with the following lemma:

\begin{lem}\label{lem:db-2}{ $f : X \to Y$ be the double cover of  $Y$ determined by 
$(B, {\mathcal L})$ as in Lemma~\ref{lem:db-1}. Let $Z$ be a smooth subvariety of $Y$ such that $(i)$ 
$\dim Z > 0$ and $(ii)$ $Z \not\subset B$. We denote the inclusion morphism
$Z \hookrightarrow Y$ by $\iota$. If there exists a divisor $B_1$ on $Z$ such that
\begin{itemize}
 \item $\iota^*B = 2B_1$ and
 \item $\iota^*{\mathcal L} \sim B_1$,
\end{itemize}
then the preimage $f^{-1}Z$ splits into two irreducible components $Z^+$ and $Z^-$.
}
\end{lem}

\proof Let $f|_{f^{-1}(Z)} : f^{-1}(Z) \to Z$ be the induced morphism. $f^{-1}(Z)$ is realized 
as a hypersurface in the total space of $\iota^*L$ as in usual manner (see \cite[Chapter I, \S 17]{bpv}, for example).
Our condition implies that $f^{-1}(Z)$ is reducible. Since $\deg f = 2$, our statement holds.
\proofend

\begin{lem}\label{lem:db-3}{Let $Y$ be a smooth projective variety, let 
$\sigma : Y \to Y$ be an involution on $Y$, let $R$ be a smooth irreducible
divisor on $Y$ such that $\sigma|_R$ is the identity, and let $B$ be a reduced divisor on $Y$
such that $\sigma^*B$ and $B$ have no common component.

If there exists a $\sigma$-invariant divisor $D$ on $Y$ (i.e., $\sigma^*D = D$) such that

\begin{itemize}
\item $B+D$ is $2$-divisible in $\Pic(Y)$, and
\item $R$ is  not contained in $\Supp(D)$,
\end{itemize}
then there exists a double cover $f: X \to Y$ branched at $2(B + \sigma^*B)$ such that
$R$ is a splitting divisor with respect to $f$.

Moreover, if there is no $2$-torsion in $\Pic(Y)$, then $R$ is a splitting divisor with 
respect to $B + \sigma^*B$.
}
\end{lem}

\proof By our assumption and $Y$ is projective, there exists a divisor $D_o$ on $Y$ such that 
\begin{enumerate}
\item
$R$ is not contained in $\Supp(D_o)$, and

\item $B+D \sim 2D_o$.
\end{enumerate}
Hence $B+\sigma^*B \sim 2(D_o + \sigma^*D_o - D)$ Let $f : X \to Y$ be a double cover determined
by $(Y, B + \sigma^*B, D_o + \sigma^*D_o - D)$ and let $\iota :
R \hookrightarrow Y$ denote the inclusion morphism. Since $\sigma |_R = \id_R$,
\[
\iota^*B = \iota^*\sigma^*B, \qquad \iota^*(D_o - D) = \iota^*(\sigma^*D_o - D),
\]
we have
\begin{eqnarray*}
\iota^*B & \sim & \iota^*(2D_o - D) \\
  & = & \iota^*D_o + \iota^*(\sigma^*D_o - D) \\
  & = & \iota^*(D_o + \sigma^*D_o-D).
  \end{eqnarray*}

Hence, by Lemma~\ref{lem:db-2}, $R$ is a splitting divisor with respect to $f$. Moreover,
if there is no $2$-torsion in $\Pic(Y)$, $f$ is determined by $B+\sigma^*B$. Hence $R$ is 
a splitting divisor with respect to $B + \sigma^*B$.
\proofend

\begin{prop}\label{prop:main}
{
Let $p_d : S_d \to \Sigma_d$ and $q_d : W_d \to \Sigma_d$ be the double covers as in
Introduction. If there exists a $\sigma_{p_d}$-invariant divisor $D$ on $S_d$ such that
$s^+ + D$ is $2$-divisible in $\Pic(S_d)$, $\sigma_{p_d}$ being the covering transformation of
$p_d$, then $T_d$ is a splitting divisor with respect to $\Delta_{0,d} + \Delta$.
}
\end{prop}

\proof Let $\psi_1$ and $\psi_2$ be rational function on $\Sigma_d$ such that
$\CC(W_d) = \CC(\Sigma_d)(\sqrt{\psi_1})$ and
 $\CC(S'_d) (= \CC(S_d)) = \CC(\Sigma_d)(\sqrt{\psi_2})$, respectively. Note that
 $(\psi_1) = \Delta_{0,d} + \Delta + 2D_1$ and $(\psi_2) = \Delta_{0,d} + T_d + 2D_2$ for 
 some divisors $D_1$ and $D_2$ on $\Sigma_d$. Let $X'_d$ be the 
 $\CC(\Sigma_d)(\sqrt{\psi_1}, \sqrt{\psi_2})$-normalization of $\Sigma_d$ and let $\tilde {p}_d: X_d
 \to S_d$ be
 the induced double cover of $S_d$ by the quadratic extension
 $\CC(\Sigma_d)(\sqrt{\psi_1}, \sqrt{\psi_2})/\CC(\Sigma_d)(\sqrt{\psi_2})$ and let 
 $\tilde{\mu} : X_d \to X'_d$ be the induced morphsim. $X'_d$ is  a bi-double
 cover of $\Sigma_d$ as well as a double cover of both $W_d$ and $S'_d$. We denote the induced
 covering morphisms by $\tilde{q}_d : X'_d \to W_d$.
 \[
 \begin{CD}
	W_d@<{\tilde{q}_d}<<X'_d@<{\tilde{\mu}}<<X_d\\
	@V{q_d}VV @V{\tilde {p}'_d}VV@VV{\tilde{p}_d}V \\
	\Sigma_d@<p'_d<<S_d@<{\mu}<<S_d
\end{CD}
\]

 Since 
 \[
 (p_d^*\psi_1) = 2O + s^+ + s^- + 2p_d^*D_1, p_d = p'_d\circ\mu 
 \]
 and 
 \[
 (q_d^*\psi_2) = 2\Delta_{0, d/2} + q_d^*T_d + 2q_d^*D_2,
 \]
 the branch loci of $\tilde {p}_d$ and $\tilde{q}_d$ are $p_d^*\Delta = s^+ + s^-$ and
 $q_d^*T_d$, respectively. 
 Put
 \[
 R:= (p_d^*T_d)_{red}\setminus (\mbox{the exceptinonal set of $S_d \to S'_d$}).
 \]
 Let $(T_d)_{sm}$ be the smooth part of $T_d$. Since $(q_d\circ\tilde {q}_d\circ\tilde{\mu})^*(T_d)_{sm}
 = (p_d\circ\tilde{p}_d)^*(T_d)_{sm}$, 
  one can check that $T_d$ is
 a splitting curve with respect to $\Delta_{0,d} + \Delta$ if and only if 
 $R$ is a splitting curve with respect to $s^+ + s^-$. Now by Lemma~\ref{lem:db-3}, our
 statement follows. \proofend

We are now in position to prove Theorem~\ref{thm:qr-1}

 \medskip
 
  \textsl{Proof of Theorem~\ref{thm:qr-1}} We first note that the algebraic equivalence 
  $\approx$ and the linear equivalence $\sim$ coincides on $S_d$ by Lemma~\ref{lem:non-etale}. 
  
  \textbf{The case of $g \ge 2$.}
   Let $s_0$ be an element in $\MW({\mathcal J}_{S_d})$ such that
   $2s_0 = s^+$ on $\MW({\mathcal J}_{S_d})$. By \cite{shioda99}, there exists a divisor $D$ on $S_d$
   which gives $s_0$. By \cite{shioda99}, $D$ satisfies  a relation
   \[
   2D \sim s^+ + (2DF - 1)O + \alpha {\mathfrak f} + \Xi,
   \]
   where $\Xi$ is a divisor whose irreducible components consist of those of singular fibers not
   meeting $O$. By our assumption on the singularity of $T_d$, we can infer that any irreducible
   component of $\Xi$ is $\sigma_{p_d}$-invariant.
  As $\sigma_{p_d}^*O = O$, $\sigma_{p_d}^*{\mathfrak 
   f} = {\mathfrak f}$, by
   Proposition~\ref{prop:main}, our statement follows.

   \textbf{The case  of $g = 1$.}    Let $s_0$ be  an element in $\MW({\mathcal J}_{S_d})$ such that
    $2s_0 = s^+$. By Theorem~\ref{thm:shioda} and Corollary~\ref{cor:mw}, we have
   \[
   2s_0 -  s^+  \in Tr(\varphi).
   \]
   Let $\phi: \MW({\mathcal J}_{S_d}) \to \NS_{\QQ}(:= \NS(S_d)\otimes\QQ)$ be the homomorphism given 
   in \cite[Lemmas 8.1 and 8.2]{shioda90}. Note that there will be no harm in considering $\NS_{\QQ}$
   since $\NS(S_d)$ is torsion free.  By \cite[Lemmas 8.1 and 8.2]{shioda90}, $\phi(s)$ satisfies the
   following properties:
   
   \begin{enumerate}
   
   \item[(i)] $\phi(s) \equiv s \bmod Tr(\varphi)_{\QQ}(:= Tr(\varphi)\otimes\QQ)$.
   
   \item[(ii)] $\phi(s)$ is orthogonal to $Tr(\varphi)$.
   
   \end{enumerate}
   
    Explicitly $\phi(s)$ is  given by
   \[
   \phi(s) = s - O - (sO + \chi({\mathcal O}_{S_d})){\mathfrak f} + \mbox{the correction terms}, 
   \]
   where ${\mathfrak f}$ denote the fiber of $\varphi$. The correction terms is a $\QQ$-divisor
   arising from reducible singular fiber in the following way:
   
   Let ${\mathfrak f}_v$ be a singular fiber over $v \in \PP^1$ and let $\Theta_{v,0}$ be
   the irreducible component with $O\Theta_0 = 1$.
   
   \begin{itemize}
   
   \item If $s$ meets $\Theta_{v,0}$, then 
   there is no correction term from ${\mathfrak f}_v$.
   
   \item If $s$ does not meet $\Theta_{v,0}$, the correction term from ${\mathfrak f}_v$ is
   as follows:
   
   Let $\Theta_{v,1}, \ldots, \Theta_{v,r_v -1}$ denote irreducible components of ${\mathfrak f}_v$ other
   than $\Theta_{v,0}$ and let $A := ( (\Theta_{v,i} \Theta_{v,j}))$ be the intersection matrix of
   $\Theta_{v,1}, \ldots, \Theta_{r_v -1}$.  With these notation,  the correction term is
\[
\sum_{i}(\Theta_{v,1},\dots,\Theta_{v, r_v-1})(-A^{-1})\left (\begin{array}{c} 
 s\Theta_{v,1} \\
\cdot \\
s\Theta_{v, r_v-1}
\end{array} \right ).
\]

\end{itemize}   
  By our assumption, 
\[
\phi(s^+) = s^+ - O - \chi({\mathcal O}_{S_d}){\mathfrak f}.
\]
Put
\[
\phi(s_0) = s_0 - O - (s_oO + \chi({\mathcal O}_{S_d})){\mathfrak f} + \sum_{v \in \Red} \mbox{Corr}_v,
\]
where $\Red = \{ v\in \PP^1 | \mbox{$\varphi^{-1}(v)$ is reducible.}\}$ and $\mbox{Corr}_v$ 
denotes the correction term arising from the singular fiber ${\mathfrak f}_v$.  Since $2s_0 - s^+ \in Tr(\varphi)$,
$\phi(2s_0) - \phi(s^+) = 0$. Hence
\[
(\ast) \quad 2s_0 - s^+ \sim_{\QQ} O + (2s_0O + \chi({\mathcal O}_{S_d})){\mathfrak f} +
2 \sum_{v \in \Red} \mbox{Corr}_v.
\]
Thus 
\[
 2\sum_{v \in \Red} \mbox{Corr}_v \sim_{\QQ} E,
\]
for some element  $E \in Tr(\varphi)$.  

\medskip

\textbf{Claim. } $2\sum_{v \in \Red} \mbox{Corr}_v \in Tr(\varphi)$.

\medskip

\textsl{Proof of Claim.} 
We first note that $2\sum_{v \in \Red} \mbox{Corr}_v = E$ in
$Tr(\varphi)_{\QQ}$. 
Since $O, {\mathfrak f}$ and all the irreducible
components of reducible singular fibers which do not meet $O$ form a basis of 
the free $\ZZ$-module $Tr(\varphi)$ as well as the $\QQ$-vector space $Tr(\varphi)_{\QQ}$,
$E$ is expressed as a $\ZZ$-linear combination of these divisors. As $\mbox{Corr}_v$ is
a $\QQ$-linear combination of the irreducible
components of reducible singular fibers which do not meet $O$,
if  $2\sum_{v \in \Red} \mbox{Corr}_v \not\in Tr(\varphi)$, then  
we have a nontrivial relation among $O, {\mathfrak f}$ and all the irreducible
components of reducible singular fibers which do not meet $O$.  This leads us to a contradiction. 
\proofend

\medskip

By Claim, we have
\begin{enumerate}

\item[(i)] $\mbox{Corr}_v = 0$ if the singular fiber over $v$ is 
of type either $I_n$ ($n$: odd), $IV$ or $IV^*$ and

\item[(ii)]  if $\mbox{Corr} \neq 0$,  one can write $\mbox{Corr}_v$ in
such a way that
\[
\mbox{Corr}_v = \frac 12 D_{1, v} + D_{2, v}, 
\]
where $D_{1,v}, D_{2,v} \in Tr(\varphi)$ and $D_{1,v}$ is reduced.

\end{enumerate}

  Since $s_0 + \sigma_{p_d}^*s_0 \in Tr(\varphi)$, we have
\[
\frac 12 (D_1 + \sigma_{p_d}^*D_1) \in Tr(\varphi).
\]
Therefore we infer that we can rewrite $D_1$ in such a way that
\[
D_1 = D'_1 + \sigma_{p_d}^*D'_1 + D''_1,
\]
where 
\begin{itemize}

\item $D'_1 \neq D'_1$ and there is no common component between $D'_1$ and $\sigma_{p_d}^*D'_1$,
 and
 
 \item each irreducible component of $D''_1$ is $\sigma_{p_d}$-invariant.
 
 \end{itemize}
 
 In particular, $D_1$ is $\sigma_{p_d}$-invariant. Now put
 \begin{eqnarray*}
 D &:= & O + (2s_0O + \chi({\mathcal O}_{S_d}) - 2\left [(2s_0O + \chi({\mathcal O}_{S_d}))/2\right ]){\mathfrak f} + D_1 \\
 D_o & := & s_0 - \left [(2s_0O + \chi({\mathcal O}_{S_d}))/2\right ]{\mathfrak f} - D_2.
 \end{eqnarray*}
 Then the relation $(\ast)$ becomes
 \[
 s^+ + D \sim 2D_o.
 \]
   As $\sigma_{p_d}^*O = O$, $\sigma_{p_d}^*{\mathfrak f} = {\mathfrak f}$, by
   Proposition~\ref{prop:main}, our statement follows.

   \medskip

    We now go on to prove the converse. Choose affine open subsets
   $V \subset W_d (= \Sigma_{d/2})$, and $U \subset \Sigma_{d}$ as follows:
   
   \begin{enumerate}
   
   \item[(i)] Both $U$ and $V$ are $\CC^2$.
   
   \item[(ii)] Let $(t,x)$ and $(\tilde{t}, \zeta)$  be affine coordinates of $U$ and $V$, respectively. Then
   $q_d$ is given by
   \[
   q_d : (\tilde{t}, \zeta) \mapsto (t, x) = (\tilde{t}, \zeta^2 + f(t)),
   \]
   where $f(t)$ is a polynomial of degree $\le d$. With respect to the coordinate $(t, x)$, 
   $\Delta_{q_d}$ is given by $\{x = \infty\} \cup \{x - f(t) = 0\}$
   
   \end{enumerate}
   
   Since  $T_d$ is a splitting curve with respect to 
   $q_{d} : W_d \to \Sigma_{d}$,  we have $q_d^*T_d = T^+ + T^-$. As $T^{\pm} \sim
   (2g+1)(\Delta_{0, d/2} + d F_{d/2})$,  we may assume that $T^{\pm}$  are given by the 
   following equations on $V$:
   \begin{eqnarray*}
   T^+ : F(x, t) + \zeta G(x,t) & = & 0 \\
   T^- : F(x, t) - \zeta G(x, t) & =  & 0,
   \end{eqnarray*}
   where
   \[
   F(x, t) = \sum_{l = 0}^g a_{2l+1}(t) (x - f(t))^{g-l}, 
   G(x, t) =  \sum_{l = 0}^g a_{2l}(t) (x - f(t))^{g-l}, a_0 = 1,
   \]
   \begin{itemize}
   
  \item  $x = \zeta^2 + f(t)$, and
  
  \item $a_{2l}(t)\,  (l =1, \ldots, g)$ and $a_{2l+1}(t) \, (l = 0, \ldots g)$ are polynomials 
   with $\deg a_{2l}(t) \le dl$ and $\deg a_{2l+1}(t) \le d(2l+1)/2$, respectively . 
   
   \end{itemize}This implies that $T_d$ is given by
   a defining equation of the form
   \[
   F(x, t)^2  - (x - f(t))G(x, t)^2 = 0.
   \]
   
   On the other hand , the generic fiber of $\varphi_d : S_d \to \PP^1$ is given by
   \[
   y^2 = F(x, t)^2  - (x - f(t))G(x, t)^2.
   \]
   By considering the divisors of the rational functions on the generic fiber $(S_d)_{\eta}$ given
   by $y - F(x, t)$ an d$y + F(x, t)$ and the right hand side of the above equation, we infer that
     the sections given by $(f(t), \pm a_{2g+1}(t))$ is $2$-divisible 
   in  $MW({\mathcal J}_{S_d})$. AS $s^{\pm}$ are nothing but these sections, our statement follows.

 

\section{Proof of Theorem~\ref{thm:application}}

We first show that $2$-divisibility of $s^+$ in $\MW({\mathcal J}_{S_d})$ follows from
the existence of a $D_{2n}$-cover for one odd number $n$. 

Suppose that there exists a $D_{2n}$-cover $\pi_d : {\mathcal X}_d \to \Sigma_d$ branched at
$2(\Delta_{0, d} + \Delta) + nT_d$ for some $n$. Let $\beta_1(\pi_d) : D({\mathcal X}_d/\Sigma_d) \to \Sigma_d$ be
the double cover canonically determined by $\pi_d : {\mathcal X}_d \to \Sigma_d$. As
the branch locus of $\beta_1(\pi_d)$ is $\Delta_{0, d} + \Delta$,
$D({\mathcal X}_d/\Sigma_d) = W_d$ and $\beta_1(\pi_d) = q_d$. By Corollary~\ref{cor:di-nec}, $T_d$
is a splitting curve with respect to $q_d$. By Theorem~\ref{thm:qr-1}, 
$s^+$ is $2$-divisible in $\MW({\mathcal J}_{S_d})$.

Conversely, suppose that $s^+$ is $2$-divisible in $\MW({\mathcal J}_{S_d})$. By Theorem~\ref{thm:qr-1},
$T_d$ is a splitting curve with respect to $q_d$. Hence we infer that $q_d^*T_d$ is of the
form $T^+ +T^-$. Put
\[
T^+ \sim a\Delta_{0, d/2} + bF_{d/2}.
\]
Since
\[
q_d^*T_d \sim (2g+1)(2\Delta_{0,d/2} + dF_{d/2},\,\, \sigma_{q_d}^*T^+ = T^-\,\, \sigma_{q_d}^*\Delta_{0,d/2}
= \Delta_{0,d/2}\,\, \mbox{and}\,\, \sigma_{q_d}^*F = F,
\]
 we have
\[
T^+ \sim T^- \sim (2g+1)\Delta_{0,d/2} + \frac {(2g+1)d}2 F_{d/2}.
\]
Hence by Corollary~\ref{cor:di-suf}, There exists a $D_{2n}$-cover
branched at $2\Delta + nT_d$ for any odd $n$. 


\section{Examples for the case of $g = 1$}

In this section, we consider the case of $g=1$. Namely, $\varphi: S_d \to \PP^1$ is an 
elliptic fibration over $\PP^1$ with section $O$. In this case,
the involution induced by the covering transformation coincides with the one induced by
the inversion morphism with respect to the group law on the generic fiber, $O$ being
the zero element.

Our main references are \cite{kodaira}, \cite{miranda-persson} and \cite{shioda-usui}.
As for the notation  of singular fibers, we follow Kodaira's notation (\cite{kodaira}).

Let us start with the case when $d=2$, i.e., $\varphi : S_2 \to \PP^1$ is a rational elliptic 
surface.

\begin{exmple}\label{eg:No7}{\rm (\cite[Example, p.198]{shioda-usui}) Let 
$\varphi : S_2 \to \PP^1$ be the rational elliptic surface given by the following
Weierstrass equation:
\[
y^2 = x^3 + (271350 - 98t)x^2  + t(t-5825)(t-2025)x + 36t^2(t-2025)^2,
\]
$t$ being an inhomogeneous coordinate of $\PP^1$.  $\varphi : S_2 \to \PP^1$ satisfies
the following properties:
\begin{enumerate}

\item[(i)] $\varphi$ has $3$ singular fibers over
$t = 0, 2025, \infty$, of which types are  of type $I_2$ over $t = 0, 2025$ and type
$III$ over $t = \infty$.

\item[(ii)] $\MW({{\mathcal J}_{S_2}})$ has no torsion.

\end{enumerate}

We infer that $T_2$ on $\Sigma_2$ given by the right hand side of the 
Weierstrass equation  has $3$ nodes (see \cite[Table 6.2, p.551]{miranda-persson}) from
(i) and is irreducible from (ii).  In order to give an example, we 
use three sections given by \cite{shioda-usui} as follows:
\[
s_0: (0, 6t^2 - 12150t), \,\, s_1: (-32t, 2t^2 - 6930t), \,\, s_2: (-20t, 4t^2 - 4500t).
\]
Let $\langle\, , \, \rangle$ be the height pairing defined in \cite{shioda90}.  Then we have
\[
\langle  s_0, s_0 \rangle = \frac 12, \,\, \langle s_i, s_i \rangle = 1\, (i = 1,2),\,\,
\langle s_1, s_2 \rangle = 0, 
\]
and there is no other section $s$ with $\langle s, s \rangle = 1/2$ other than $\pm s_0$.

The sections given by $2s_0$ and $s_1+s_2$ are
\[
2s_0 =\left (\frac 1{144}t^2  + \frac {1231}{72}t - \frac{5143775}{144}, 
        -\frac{1}{1728}t^3- \frac{2335}{576}t^2 + \frac{13493375}{576}t -
       \frac{29962489375}{1728} \right )
\]  
\[
s_1+s_2=\left (\frac1{36}t^2 + \frac{435}2 t - \frac{921375}4, 
 -\frac 1{216}t^3 - \frac{1181}{24}t^2 - \frac{41625}8 t + \frac{373156875}8\right )
 \]
 Since $s_1 + s_2 \neq -2s_0$, we infer that $2s_0$ is $2$-divisible, while $s_1 + s_2$ is
 not $2$-divisible. Also, both $2s_0$ and $s_1 + s_2$ do not meet the zero section $O$. Let $\Delta^{(1)}$ and $\Delta^{(2)}$ be the sections which are the images of 
 $2s_0$ and $s_1 + s_2$ in $\Sigma_2$, respectively.
 Put
 \[
 B_1 = \Delta_{0, 2} + \Delta^{(1)} + T_2 , \quad B_2 = \Delta_{0,2} + \Delta^{(2)} + T_2.
 \]
 Oen can check that, for each $i$, $\Delta^{(i)}$ and $T_2$ meet $3$ distinct smooth points of 
 $T_2$ in such a way
 that the  intersection multiplicity at each point is $2$. Let $q_2^{(i)} : W^{(i)}_2 \to \Sigma_2$ 
 $(i = 1, 2)$ be
 the double covers with branch locus $\Delta_{0, 2} + \Delta^{(i)}$ $(i = 1, 2)$, respectively. Then
 $T_2$ is a splitting curve with respect to $q_2^{(1)}$, but not with respect to $q_2^{(1)}$.
 Hence by Corollary~\ref{cor:zpair}  there exists no homeomorphism $h: \Sigma_2 \to \Sigma_2$ such
 that $h(B_1) = B_2$.
}
\end{exmple}

\begin{exmple}\label{eg:No10}{\rm (\cite[Example, p. 210]{shioda-usui}) Let 
$\varphi : S_2 \to \PP^1$ be the rational elliptic surface given by the following
Weierstrass equation:
\[
y^2 = x^3 + (25t+9)x^2 + (144t^2 + t^3)x + 16t^4,
\]
$t$ being an inhomogeneous coordinate of $\PP^1$.  (Note that
the original Weierstrass equation in \cite{shioda-usui} is
$y^2 -6xy= x^3 + 25tx^2  + (144t^2 + t^3)x + 16t^4$. We change the equation slightly.)
$\varphi : S_2 \to \PP^1$ satisfies
the following properties:
\begin{enumerate}

\item[(i)] $\varphi$ has $2$ singular fibers over
$t = 0, \infty$, of which types are  of type $I_4$ over $t = 0$ and type
$III$ over $t = \infty$.

\item[(ii)] $\MW({\mathcal J}_{S_2})$ has no torsion.

\end{enumerate}

Likewise Example~\ref{eg:No7}, we infer that $T_2$ on $\Sigma_2$ has one $a_3$-singularity and one $a_1$-singularity and is irreducible. 
In order to give another example, we 
use three sections given by \cite{shioda-usui} as follows:
\[
s_0: (0, 4t^2), \,\, s_1: (-16t, -48t), \,\, s_2: (-15t, t^2+45t).
\]
Let $\langle\,\,  , \, \rangle$ be the height pairing defined in \cite{shioda90}.  Then we have
\[
\langle  s_0, s_0 \rangle = \frac 12, \,\, \langle s_i, s_i \rangle = \frac 34\, (i = 1,2),\,\,
\langle s_1, s_2 \rangle = -\frac 14, 
\]
and there is no other section $s$ with $\langle s, s \rangle = 1/2$ other than $\pm s_0$.
The sections given by $2s_0$ and $s_1+s_2$ are
\[
2s_0 =\left (\frac 1{64}t^2 - \frac{41}2 t + 315, -\frac{1}{512}t^3 - \frac{55}{32}t^2 + \frac{2637}{8} - 5670 \right )
\]  
\[
s_1+ s_2=\left (t^2 + 192t + 8640, -t^3-301t^2 - 27936t-803520 \right )
 \]
Since $s_1 + s_2 \neq -2s_0$, we infer that $2s_0$ is $2$-divisible, while $s_1 + s_2$ is
 not $2$-divisible. Also, both $2s_0$ and $s_1 + s_2$ do not meet the zero section $O$. Let $\Delta^{(1)}$ and $\Delta^{(2)}$ be the sections which are the images of 
 $2s_0$ and $s_1 + s_2$ in $\Sigma_2$, respectively.
 Put
 \[
 B_1 = \Delta_{0, 2} + \Delta^{(1)} + T_2 , \quad B_2 = \Delta_{0,2} + \Delta^{(2)} + T_2.
 \]
For the pair $(B_1, B_2)$,  similar properties to those in Example~\ref{eg:No7} hold.
}
\end{exmple}

\begin{rem}{\rm Note that the right hand sides of Weierstrass equations of Examples~\ref{eg:No7} and ~\ref{eg:No10} can be rewritten as follows:
\begin{eqnarray*}
 && x^3 + (271350 - 98t)x^2  + t(t-5825)(t-2025)x + 36t^2(t-2025)^2 \\
 &= & \left \{ \left( \frac 1{12} t - \frac{5825}{12}\right ) x + (6t^2 - 12150t)\right \}^2 +
 x^2\left ( x - \frac 1{144} t^2 - \frac{1231}{72} t  + \frac {5143775}{144} \right ) , \\
 & & x^3 + (25t+9)x^2 + (144t^2 + t^3)x + 16t^4 \\
& = & \frac 1{64}\left \{ (t + 144)x+ 32t^2\right \}^2 + x^2 \left (x - \frac1{64}t^2 + 
\frac {41}2 t - 315 \right )
 \end{eqnarray*}
Since the sections of $\Sigma_2$ given by 
$x - \frac 1{144} t^2 - \frac{1231}{72} t  + \frac {5143775}{144} = 0$ and 
$x - \frac1{64}t^2 + 
\frac {41}2 t - 315= 0$ are $\Delta^{(1)}$'s in Examples~\ref{eg:No7} and ~\ref{eg:No10}, respectively, we see that $\Delta^{(1)}$'s in both of examples are splitting curves
with respect to $q_2^{(1)}$.
 }
\end{rem}

We show that an infinitely many examples of pairs $(B_1, B_2)$ with similar properties to
those in Examples~\ref{eg:No7} and \ref{eg:No10}. Let $\varphi : {\mathcal E} \to \PP^1$ be an
elliptic surface with section $O$ given by 
\[
(\ast) \quad y^2 = x^3 + a_2(t)x^2 + a_4(t)x + a_6(t), 
\]
where $a_2(t), a_4(t), a_6(t) \in \CC(t)$. We assume that $\varphi : {\mathcal E} \to \PP^1$
satisfies the following two properties:

\begin{itemize}

\item $\MW({\mathcal J}_{\mathcal E})$ has no torsion.

\item There exists sections $s_1$ and $s_2$ in $\MW({\mathcal J}_{\mathcal E})$ such that $s_1$ is $2$-divisible
$\MW({\mathcal J}_{\mathcal E})$, while $s_2$ is not $2$-divisible in $\MW({\mathcal J}_{\mathcal E})$.

\end{itemize}

Let $\nu : C \to \PP^1$ be a surjecticve morphism from a smooth curve $C$ to $\PP^1$ and let $\varphi_{\nu} : {\mathcal E}_{\nu}
\to \PP^1$ be the induced elliptic fibration obtained by the pull back by $\nu$. 

\begin{lem}\label{lem:pullback}{If we choose $C$ and $\nu$ generally enough, the following properties are
satisfied:
\begin{enumerate}

\item $\MW({\mathcal J}_{{\mathcal E}_{\nu}})$ has no  torsion of oder $n$.

\item Let $s_{\nu, 1}$ and $s_{\nu, 2}$ be the sections in 
$\MW({\mathcal J}_{{\mathcal E}_{\nu}})$ induced
by $s_1$ and $s_2$, respectively.   $s_{\nu,1}$ is $2$-divisible
$\MW({\mathcal J}_{{\mathcal E}_{\nu}})$, while $s_{\nu, 2}$ is not $2$-divisible in $\MW({\mathcal J}_{{\mathcal E}_{\nu}})$.

\end{enumerate}

}
\end{lem}

\proof  We first note that $\MW({\mathcal J}_{\mathcal E})$ and  $\MW({\mathcal J}_{{\mathcal E}_{\nu}})$ are regarded
as the set of $\CC(t)$- and $\CC(C)$- rational points of the elliptic curve $(\ast)$, respectively. 

 1. Let $K_n$ be the field extension of $\CC(t)$ obtained 
by addding all $n$-torsion points.  
If we choose $C$ and $\nu$ generically enough so that $\CC(C) \cap K_n = \CC(t)$, then
$\MW({\mathcal J}_{{\mathcal E}_{\nu}})$ has no $n$-torsion.

2. The duplication formula for the group law of an elliptic curve shows that
$s_2$ becomes $2$-divisible over a certain field extension  $K$ of $\CC(t)$ of degree at most $4$. 
Hence if choose $C$ and $\nu$ generically enough so that $\CC(C) \cap K = \CC(t)$, then
$s_{\nu, 2}$ is not $2$-divisible in 
$\MW({\mathcal J}_{{\mathcal E}_{\nu}})$.
\proofend 

\medskip

We now give examples of pairs $(B_1, B_2)$ on $\Sigma_{2n}$ as in Examples~\ref{eg:No7} and
~\ref{eg:No10} 
for any $n$.

\begin{exmple}\label{eg:infinite}{\rm Let $\varphi : S_2 \to \PP^1$ be the rational elliptic surface
in Example~\ref{eg:No7} or \ref{eg:No10}. Let $\nu_n : \PP^1 \to \PP^1$ be a generic
rational map of degree $n$. Let $\varphi_{\nu_n} : {\mathcal E}_{\nu_n} \to \PP^1$ be the induced
relatively minimal elliptic fibration obtained via the pull-back of $\nu_n$ from $\varphi$. As  we have
seen at the beginning of this section, ${\mathcal E}_{\nu_n}$ is the canonical resolution of a double
cover $p'_{2n}: {\mathcal E}'_{\nu_n} \to \Sigma_{2n}$ of $\Sigma_{2n}$ such that
\begin{itemize}
\item the branch locus of $p'_{2n}$ is of the form $\Delta_{0,2n} + T_{2n}^{(\nu)}$, 
$T_{2n}^{(\nu)} \sim 3(\Delta_{0,2n} + 2nF_{2n})$, and $T_{2n}^{(\nu)}$ has at worst
simple singularities, and

\item the involution induced by the covering transformation of $p'_{2n}$ coincides with
the one induced by the inversion morphsim of the group law on the generic fiber.
 
\end{itemize}
By Lemma~\ref{lem:pullback}, we may assume that $\MW({\mathcal J}_{{\mathcal E}_{\nu_n}})$ has
no $2$-torsion. This implies that $T_{2n}^{(\nu_n)}$ is irreducible. 
Put
\begin{eqnarray*}
s_{1, \nu_n} & := & \mbox{the section arising from $2s_o$} \\
s_{2, \nu_n} & := & \mbox{the section arising from $s_1 + s_2$.}
\end{eqnarray*}
and let $\Delta_{\nu_n}^{(1)}$ and $\Delta_{\nu_n}^{(1)}$ be the sections of $\Sigma_{2n}$ induced by
$s_{1, \nu_n}$ and $s_{2, \nu_n}$, respectively. Put
 Put
 \[
 B_1^{(\nu_n)} = \Delta_{0, 2} + \Delta_{\nu_n}^{(1)} + T_{2n}^{(\nu_n)} , \quad B_2^{(\nu_n)} = \Delta_{0,2} + \Delta_{\nu_n}^{(2)} + T_{2n}^{(\nu_n)}.
 \]
For the pair $(B_1^{(\nu_n)}, B_2^{(\nu_n)})$,  similar properties to those in Example~\ref{eg:No7} hold.

}
\end{exmple}


\section{Application to the study of Zariski pairs}

In this section, we give two examples of Zariski pairs of sextic curves by using Corollary~\ref{cor:zpair}.

Let $C$ be a smooth conic and let $Q$ be an irreducible quartic such that

\begin{enumerate}

 \item[(i)] $Q$ has either $2a_1$ singularities or one $a_3$ singularity,
 
 \item[(ii)] $C$ is tangent to $Q$ at $4$ distinct smooth points of $Q$, and
 
 \item[(iii)] one of the $4$ intersection points of $C$ and $Q$ is an inflection point, $x$,
 of $Q$.
 
 \end{enumerate}
 
 Likewise \cite[\S 4]{tokunaga94}, we consider a rational elliptic surface related to $C+Q$ as follows:
 
 \medskip
 
 Let $\nu_1 : \PP^2_x \to \PP^2$ be a blowing-up at $x$. Let $Q_x$ and $E$ denote the
 proper transform of $Q_x$ and the exceptional divisor of $\nu$, respectively. Next let 
 $\nu_2 : \widehat{\PP}^2 \to \PP^2_x$ be a blowing-up at $Q_x\cap E$, and let 
 $\widehat{Q}, E_1$, and $E_2$ be the proper transforms of $Q_x$ and $E$ and the exceptional
 divisor of $\nu_2$, respectively.
 
 Let $\overline{l}_x$ be the proper transform of the tangent line at $x$. $\overline{l}_x$ is
 the exceptional curve of the first kind. By blowing down $\overline{l}_x$, we obtain
 $\Sigma_2$, and $\widehat{Q}$ and $C$ give rise to divisors $T_{Q}$ and $\Delta_C$, respectively,
 such that 
 $\Delta_C \sim \Delta_{0,2} + 2F_2$, $T_Q \sim 3(\Delta_{0,2} + 2F_2)$ and 
 $T_Q$ has $3a_1$ (resp. $a_3 + a_1$) singularities if $Q$ had $2a_1$ (resp. one $a_3$).
 
 Let $f' : {\mathcal E}' \to \Sigma_2$ be the double cover with branch locus
 $\Delta_{0,2} + T_Q$ and we denote its canonical resolution of ${\mathcal E}'$ by ${\mathcal E}$. 
 ${\mathcal E}$ satisfies the following properties:
 
 \begin{itemize}
 
  \item ${\mathcal E}$ has $3$ singular fibers whose configuration is
  $2I_2$ and $III$.
  
  \item The conic $C$ gives rise to two sections $s^+$ and $s^-$.
  
  \item The covering transform induced by $f'$ coincides with the one induced by
  the inversion morphism with respect to the group law.
  
 \end{itemize}
 
 Under these setting, one easily infers that
 
 \begin{prop}\label{prop:zpair}{There exists a $D_{2n}$-cover branched at
 $2C + nQ$ if and only if there exists a $D_{2n}$-cover branched at
 $2(\Delta_{0,2} + \Delta_C) + nT_Q$.
 }
 \end{prop}
 
 \begin{exmple}\label{eg:zpair}{\rm Let ${\mathcal E}$ be the  rational elliptic surface as 
 in Example~\ref{eg:No7}.  By considering the reverse process observed as above starting from
 ${\mathcal E}$, we infer that
 
 \begin{itemize}
 
 \item two sections $2s_o$ and $s_1+s_2$ in 
 Examples~\ref{eg:No7} give rise to two conics $C_1$ and $C_2$, respectively and
 
 \item the corresponding $T_2$ gives rise to an irreducible quartic $Q$ with $2a_1$ singularities.
 
 \end{itemize}

 By Proposition~\ref{prop:zpair}, we infer that $(C_1\cup Q, C_2\cup Q)$ is  a Zariski pair.

 \medskip

 We similarly obtain another Zariski pair starting from the rational elliptic surface in
 Example~\ref{eg:No10}. It is a pair of sextic curves $(C_1\cup Q', C_2\cup Q')$  such that
 
 \begin{itemize}
 \item $Q'$ is an irreducible quartic with one $a_3$ singularity.
 
 \item both of $C_i$ $(i = 1,2)$ are tangent to $Q'$ at four distinct smooth points of $Q'$.
 \end{itemize}
 }
 \end{exmple}

\begin{rem}{\rm Both of Zariski pairs in Example~\ref{eg:zpair} are in Shimada's list in \cite{shimada}.
Our reasoning why they are Zariski pairs gives another point of view different from Shimada's one. 
}
\end{rem}


\noindent Hiro-o TOKUNAGA\\
Department of Mathematics and Information Sciences\\
Graduate School of Science and Engineering\\
Tokyo Metropolitan University\\
1-1 Minami-Ohsawa, Hachiohji 192-0397 JAPAN \\
{\tt tokunaga@tmu.ac.jp}
      
  \end{document}